\documentclass[a4paper,11pt]{amsart}

\usepackage{amssymb}
\usepackage{graphicx}

\begin{document}

\title{On noise limited cellular networks}

\author{L.~Decreusefond}
\address{Institut Telecom, Telecom Paristech, CNRS LTCI\\
46, rue Barrault,  Paris- 75634,
  France}
\email{Laurent.Decreusefond@telecom-paristech.fr}

\author{P.~Martins}
\address{Institut Telecom, Telecom Paristech, CNRS LTCI\\
46, rue Barrault,  Paris- 75634,
  France}
\email{Philippe.Martins@telecom-paristech.fr}

\author{T~.T.~Vu}
\address{Institut Telecom, Telecom Paristech, CNRS LTCI\\
46, rue Barrault,  Paris- 75634,
  France}
\email{vu@telecom-paristech.fr}

\newtheorem{theorem}{Theorem}
\newtheorem{lemma}[theorem]{Lemma} \newtheorem{remark}{Remark}

\newtheorem{definition}{Definition}

\newtheorem{proposition}[theorem]{Proposition}
\newtheorem{assumption}{Assumption}

\newtheorem{collary}{Corollary}
\def\R{{\mathbf R}}
\def\d{\text{ d}}
\newcommand{\var}{\operatorname{\text{var}}}
\newcommand{\cov}{\operatorname{\text{cov}}}
\begin{abstract}
  This paper introduces a general theoretical framework to analyze
  noise limited networks. More precisely, we consider two homogenous
  Poisson point processes of base stations and users. General model of
  radio signal propagation and effect of fading are also
  considered. The main difference of our model with respect to other
  existing models is that a user connects to his best servers but not
  necessarily the closest one.  We provide general formula for the
  outage probability. We study functionals related to the SNR as well
  as the sum of these functionals over all users per cell. For the
  latter, the expectation and bounds on the variance are obtained.
\end{abstract}




\maketitle

\section{Introduction}

\subsection{Motivation}
Cellular network is a kind of radio network consisting of a number of
fixed access points known as base stations and a large number of users
(or mobiles). Each base station covers a geometrical region known as a
cell and  serve all users in this cell. Interference and
noise are two factors annoying  communications in cellular wireless
networks. Noise is unavoidable and comes from natural
sources. Interferences come from  users and
base stations. The use of recent technologies such as SDMA (spatial
division multiple access) and MIMO (multiple input multiple output)
can reduce significantly  interferences so that we can hope
that in a near future the impact of interferences  will be negligible and noise will become the only factor
harming the network. The best case is when interferences
from other cells are perfectly canceled, the  network is then said to
be in noise limited regime. In this paper, we consider and introduce a
framework to study this kind of network.

In existing  literature, base stations (BS) locations are usually modeled
as an ideal regular hexagonal lattice. In reality,  base
stations are irregularly located, especially in an urban area, and the
cell radius is not the same for each BS. In this paper, we model the
base station locations as an homogenous Poisson point process $\Pi_B$ of
 intensity $\lambda_B$. Such a model comes from stochastic
geometry. It is sufficently versatile by changing $\lambda_B$ to cover
a wide number of real situations and it is  mathematically
tractable. For an introduction to the usage of
stochastic geometry for wireless networks performances, we refer to
\cite{StochasticGeometryandRandomGraphsfortheAnalysisandDesignofWirelessNetworks}. Theory
and number of pertinent examples can be found in \cite{StochasticGeometryandWirelessNetworksVolumeI} and
\cite{StochasticGeometryandWirelessNetworksVolumeII}. For all
theoretical details, we refer the first opus.

 To model cellular network cells,
Voronoi tessellations are frequently used. It is based on the
assumption that each user is served by the closest BS. Unfortunately,
this is not always very accurate since in real life, a mobile connects
to the \emph{best} BS it can have, i.e., the BS which offers it the
best Signal over Noise Ratio. The \emph{best} BS is not always the
closest because of the fading environment. In this paper, we analyze
the impact of fading by considering that users are served by the base
station providing the best signal power. The location of users in the
plane are modeled as another homogenous Poisson point process $\Pi_M$
of intensity $\lambda_M$.

While ancient cellular networks such as GSM and GPRS provided only
voice service and low data transmission rate, recent and emergent
wireless cellular networks such as WIMAX or LTE offer higher data rate
and other services requiring high throughput such as video calls. Each
service requires a different level of signal to noise ratio (SNR). If
the SNR does not reach a required threshold due to the radio
condition, the service can not be established or may be
interrupted. Such calls are said to be in outage. The outage
probability is one of key measurement of the network performance. One
of the aims of this paper is to determine the outage probability of
noise limited network, or equivalently the distribution of SNR, which
turns out to be equivalent to determine the distribution of the smallest path loss fading. In
fact, there have been some works dealing with the outage probability
of noise limited wireless network, but almost all of them consider the
exponent path loss model. This paper provides a general formula for
outage probability taking into account a more general  model of path loss.

Once the distribution of SNR of a user is determined, the distribution
of  functionals related to SNR can be easily derived.  In some
situations, we have to study the distribution of the sum of a
functional for all users in a cell. For example, in an OFDMA noise
limited cellular system, the number of sub channels required for a
user demanding a particular service depends on its SNR. If the total
number of sub channels of all users in a cell excesses the number of
available sub channels in this cell then at least one user is
blocked. The probability of that to happen, sometimes called
infeasibility probability, contains extremely important information on
the performance of the network.
Since it is often impossible to find the explicit probability
distribution of additive functionals, we calculate the expectation,
and bounds on the variance of such random variables.

This paper is organized as follows. In Section \ref{section: Model
  scenario}, we describe the model. In Section \ref{section:
  POISSON POINT PROCESS of pathloss fading}, we show that the path loss
fading can be viewed as a Poisson point process on the real line and
we provide a general
formula for the outage probability.  In Section \ref{section:
  Capacity user}, we calculate  average  capacity of a user and of a
cell. We also compute upper and lower bounds for their variance as
closed form expressions seem untractable. Section
\ref{section: Examples} illustrates the results obtained for some
particular situations.
\section{Model} \label{section: Model scenario}

\subsection{System model} \label{subsec: Propagation model} Consider a
BS (base station) located at $y$ with transmission power $P$ and a
mobile located at $x$. The mobile's received signal has average power
$L(y-x)P$ where $L$ is the path loss function. We assume that $L$ is
measurable function on $\R^2$.
The most used path loss function is the so-called path loss exponent model
\begin{equation*}
  L(z)=K|z|^{-\gamma},
\end{equation*}
where $|z|$ refers to the Euclidean norm of $z$. This function gives
raise to nice closed formulas but is rather unrealistic: Close to the
BS, the signal is infinitely amplified. A more realistic  model is the
modified path loss model given by:
\begin{equation*}
  L(z)=K\min\{R_0^{-\gamma},|z|^{-\gamma}\}
\end{equation*}
where $R_0$ is a reference distance and $K$ a constant depending on
the environment. In addition to this deterministic
large scale effect, we consider the fading effect, which is by essence
random. The received
signal power from a BS located at $y$ to a mobile unit (MU for short)
located at $x$ is given by
\begin{equation*}
  P_{yx} = h_{y,\, x}L(y-x)P,
\end{equation*}
where $\{h_{y,\, x}\}_{x,y\in \R^2}$ are independent copies of a random
variable $H$. Most used fading random models are log-normal shadowing
and Rayleigh fading. The log-normal shadowing is such that $H$ is a
log-normal random variable and we can write $H\sim 10^{G/10}$ where
$G\sim \mathcal{N}(0,\sigma^2)$. The Rayleigh fading is such that $H$
is an exponential random variable of parameter $\mu$. We can also consider the Rayleigh-Lognormal composite
fading, in this case the fading is the product of the log-normal
shadowing factor and the Rayleigh fading factor. It is worth noting
that the log-normal shadowing usually improves the network performance
while Rayleigh fading usually degrades  performances. 

We  assume
that once  in the network, a mobile is attached to the BS that
provides it the best signal strength.  If the power received at this point
 is greater than some threshold $T$, 
we say that $x$ is covered. If $x$ is not covered by any BS
then a MU at $x$ can not establish a communication  and
thus is said to be in outage. 
  In the case of path loss exponent
model with no fading ($H=\text{constant}$), the best BS for given
mobile is always its nearest BS. 

We assume that the point process of BSs $\Pi_B = \{y_0,y_1,...\}$ is an
homogenous Poisson point process of intensity $\lambda_B$ on
$\mathbf{R}^2$ and that users are distributed in the plane as a Poisson
point process $\Pi_M =
\{x_0,x_1,...\}$ of intensity $\lambda_M$.

To avoid any technical difficulty, from now on, we make the following
assumptions:

\begin{assumption}\label{assumption} Assume that:
  \begin{enumerate}
  \item All random variables $h_{yx}$ ($x,y\in \mathbf{R}^2$) are
    independent.
  \item $H$ admits a probability density function $p_H$. 
Its complementary cumulative distributive function is denoted by $F_H$, i.e., $$F_H(\beta)=P(H\geq
    \beta)=\int_{\beta}^{\infty}p_H(t)\d t >0.$$
  \item Define $B(\beta)=\int_{R^2}F_H((L(z)\beta)^{-1})\d z$. Then, we
    have 
    $0<B(\beta)<\infty$ for all $\beta>0$.
  \end{enumerate}
\end{assumption}

\begin{figure}
  \centering
  \includegraphics[width=\textwidth]{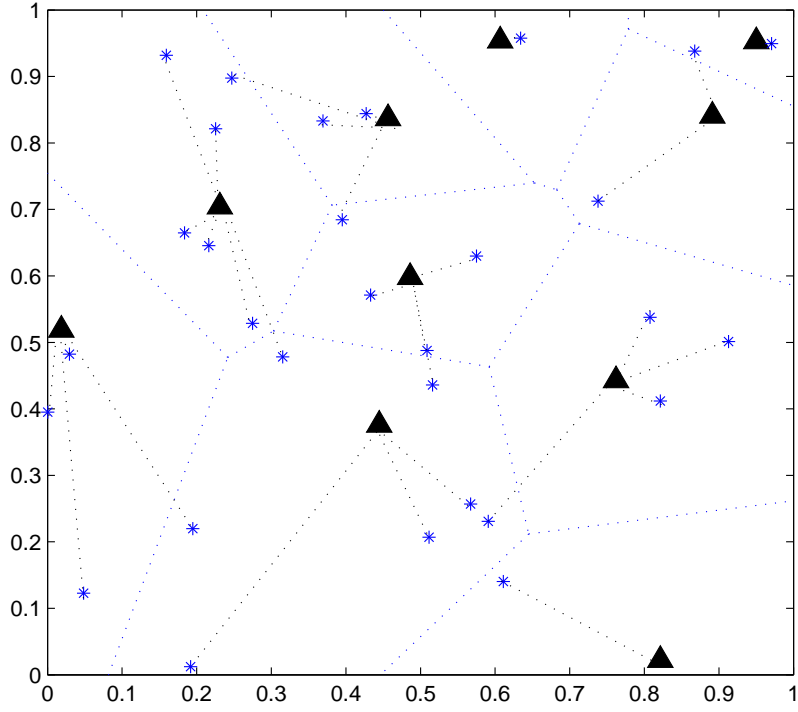}
  \caption{Triangles represent BS, plus represent MU. Dotted
    polygons are Voronoi cells induced by BS. A line between a BS and an MU means
    that the BS serves the MU. A mobile may be not served by the BS
    closest to  it, due to fading.}\label{OurModelFigure}
\end{figure}

\section{Poisson point process of path loss fading} \label{section:
  POISSON POINT PROCESS of pathloss fading}

In this section, similarly to
~\cite{AGeometricInterpretationofFadinginWirelessNetworks}, we show that the path loss fading process is a Poisson
point process on the positive half of the real line.

For each location $x$ on $\mathbf{R}^2$, consider the path loss fading
process $Sh^x = \{s_i^x\}_{i=0}^{\infty}$ where $s_i^x =
(h_{y_ix}L(y-x))^{-1}$.
\begin{proposition}\label{proposition: POISSON POINT PROCESS of path
    shadow}
  For any $x$, $Sh^x$ is a Poisson point process on $\R^{+}$
  with intensity density $d\Lambda(t)=\lambda_BB'(t)d t$. In addition,
  $B(0)=0$ and $B(\infty)=\infty$.

\end{proposition}
\begin{proof}
  Define the marked point process $\Pi_B^x =
  \{y_i,h_{y_ix}\}_{i=0}^{\infty}$. It is a Poisson point process of
  intensity $\lambda_Bd y\otimes f_H(t)d t$ because the marks are
  independent and identically distributed.  Considering the
  probability kernel $p((z,t),A) = 1_{\{(L(z)t)^{-1}\in A\}}$ for all
  Borel $A \subset \mathbf{R}^{+}$ and applying the displacement theorem
  (\cite{StochasticGeometryandWirelessNetworksVolumeI}, theorem 1.3.9),
  we obtain that the point process $Sh^x$ is a Poisson point process of
  intensity
  \begin{eqnarray*}
    \Lambda(A)=\lambda_B\int_{R^2}\int_{R}1_{(\{L(z)t)^{-1}\in
      A\}}p_H(t)\d z \d t  \cdotp
  \end{eqnarray*}
  We now show that $\Lambda([0,\beta])=\lambda_BB(\beta)$. Indeed,
  \begin{eqnarray*}
    \Lambda([0,\beta]) &=& \lambda_B\int_{R^2}\int_{R}1_{\{t\geq (\beta L(z))^{-1}\}}p_H(t)\d z \d t\\
    &=&\lambda_B\int_{R^2}F_H((\beta L(z))^{-1})\d z\\
    &=&\lambda_BB(\beta) \cdotp
  \end{eqnarray*}
  It is easy to see that $B(0)=0$ and $B(\infty)=\infty$.  Finally
  $B(\beta)$ admits the derivative:
  \begin{equation*}
    B'(\beta) = \beta^{-2}\int_{R^2}\frac{1}{L(z)}p_H((\beta
    L(z))^{-1})\d z \cdotp
  \end{equation*}
  This concludes the proof.
\end{proof}

For any point $x$, we can reorder the points of $Sh^x. $ We denote
ordered atoms of $Sh^x$ by $0\le \xi_0^x<\xi_1^x<\ldots$. The CDF and PDF of $\xi^x_m$ are easily derived
according to the property of Poisson point processes:
\begin{collary}\label{proposition CDF PDF of a_ms}
  The complementary cumulative distribution function of $\xi_m^x$ is
  given by:
  \begin{equation*}
    P(\xi_m^x>t) =
    e^{-\lambda_BB(t)}\sum_{i=0}^{m}\frac{(\lambda_BB(t))^i}{i!},
  \end{equation*}
  and its probability density function is given by
  \begin{equation}\label{equ: PEF mth si}
    p_{\xi_m^x}(t) =
    \frac{\lambda_B^{m+1}B'(t)(B(t))^m}{m!}e^{-\lambda_BB(t)} \cdotp
  \end{equation}
\end{collary}
\begin{proof}
  The event $(\xi_m^x>t)$ is equivalent to the event (in the interval
  $[0,t],$ there are at most $m$ points) and the number of points in
  this interval follows a Poisson distribution of mean
  $\lambda_BB(t)$. Thus, we have:
  \begin{eqnarray*}P(\xi_m^x>t) =
    e^{-\lambda_BB(t)}\sum_{i=0}^{m}\frac{(\lambda_BB(t))^i}{i!}\cdotp
  \end{eqnarray*}
  The PDF is thus given by
  \begin{eqnarray*}
    p_{\xi_m^x}(t) &=& -\frac{\partial}{\partial t}P(\xi_m^x>t)\\
    &=& -\lambda_BB'(t)e^{-\lambda_BB(t)} \\ &\qquad&+\sum_1^m\lambda_BB'(t)e^{-\lambda_BB(t)}\left(\frac{(\lambda_BB(t))^{i-1}}{(i-1)!}-\frac{(\lambda_BB(t))^{i}}{i!}\right)\\
    &=& \frac{\lambda_B^{m+1}B'(t)(B(t))^m}{m!}e^{-\lambda_BB(t)}
    \cdotp
  \end{eqnarray*}
The proof is thus complete.
\end{proof}

\begin{collary}\label{lemma: B(t) of path loss exponent }
  If $L(z)=K|z|^{-\gamma}$ then:
  \begin{equation*}
    B(\beta) = C.\beta^{\frac{2}{\gamma}} ,
  \end{equation*}
  where $C=\pi K^{\frac{2}{\gamma}}E(H^{\frac{2}{\gamma}})$.
\end{collary}
\begin{proof}
  The path loss function depends only on the distance from the BS to
  the user. By the change of variable $r=|z|$ and by integration by
  substitution, we have:
  \begin{eqnarray*}
    B(\beta) &=& 2\pi\int_0^{\infty}\int_0^{\infty}r1_{\{tK\beta\geq
      r^{\gamma}\}}p_H(t)\d r \d t\\
    &=& 2\pi\int_0^{\infty}p_H(t)\d t
    \int_0^{(tK\beta)^{1/\gamma}}r\d r\\
    &=&\pi(K)^{\frac{2}{\gamma}}\beta^{\frac{2}{\gamma}}\int_0^{\infty}p_H(t)t^{\frac{2}{\gamma}}\d t\\
    &=&\pi(K)^{\frac{2}{\gamma}}E(H^{\frac{2}{\gamma}})\beta^{\frac{2}{\gamma}} \cdotp
  \end{eqnarray*}
  Hence the result.
\end{proof}
Remark that this result can also be derived from
\cite{AGeometricInterpretationofFadinginWirelessNetworks}. We observe
that the distribution of the point process $Sh^x$ does depend only on
$E(H^{\frac{2}{\gamma}})$ but not on the distribution of fading $H$
itself. This phenomenon can be explained as in
\cite{InterferenceinLargeWirelessNetworks}(page 159). If the fading is
log-normal shadowing, i.e $H\sim 10^{G/10}$ where $G\sim
\mathcal{N}(0,\sigma^2)$ then
$E(H^{\frac{2}{\gamma}})=e^{\frac{2\sigma_1^2}{\gamma^2}}$ where
$\sigma_1=\frac{\ln(10)\sigma}{10}$. If the fading is Rayleigh fading,
i.e $H\sim \exp(\mu)$ then
$E(H^{\frac{2}{\gamma}})=\Gamma(\frac{2}{\gamma}+1,0)\mu^{-\frac{2}{\gamma}}$
where $\Gamma(a,b) = \int_b^{\infty}t^{a-1}e^{-t}\d t$ is the upper
incomplete gamma function.

Similarly to the distance to $m$-th nearest BS (which can be found in
\cite{DistanceDistributionsinFiniteUniformlyRandomNetworks}), the
distribution of $m$-th less strong path loss fading $\xi_m^x$ can be
characterized as follows:
\begin{collary}
  If $L(z)=K|z|^{-\gamma}$, $\xi_m^x$ is distributed according to the
  generalized Gamma distribution:
  \begin{equation*}
    p_{\xi_m^x}(t) =
    \frac{2}{\gamma}(\lambda_BC)^{m+1}t^{\frac{2}{\gamma}(m+1)}\frac{e^{-\lambda_BCt^{\frac{2}{\gamma}}}}{m!}\cdotp
  \end{equation*}
\end{collary}
\begin{proof}
  This is a consequence of Proposition \ref{proposition CDF PDF of
    a_ms} and Lemma \ref{lemma: B(t) of path loss exponent }.
\end{proof}
We can also investigate more general and realistic path loss model.
\begin{collary}
  If $L(z)=K\min\{R_0^{-\gamma},|z|^{-\gamma}\}$ then:
  \begin{equation}\label{equ: corr: modif exponent}
    B(\beta)=C_1\beta^{\frac{2}{\gamma}}\int_{\frac{
        R_0^{\gamma}}{\beta K}}^{\infty}t^{\frac{2}{\gamma}}p_H(t)\d t ,
  \end{equation}
  where $C_1=\pi K^{\frac{2}{\gamma}}$. In addition, we have:
  \begin{equation} \label{equ: corr: modif exponent path Bphay}
    B'(\beta) = \frac{2}{\gamma}\beta^{-1}B(\beta)+\pi
    R_0^2p_H\left(\frac{R_0^{\gamma}}{K\beta}\right) \cdotp
  \end{equation}
  If the fading is lognormal shadowing $H\sim 10^{G/10}$ where $G\sim
  \mathcal{N}(0,\sigma^2)$ then we have:
  \begin{equation*}
    B(\beta) =
    C_1\beta^{\frac{2}{\gamma}}e^{(\frac{2\sigma_1}{\gamma})^2}Q\left(\frac{-\ln\beta-\ln(KR_0^{-\gamma})}{\sigma_1}-\frac{2\sigma_1}{\gamma}\right),
  \end{equation*}
  where $Q(a)=\frac{1}{\sqrt{2\pi}}\int_a^{\infty}e^{-u^2/2} \d u$ is
  the \emph{Q}-function and $\sigma_1=\frac{\sigma\ln 10 }{10}$.  If
  the fading is Rayleigh $H\sim \exp(\mu)$ then
  \begin{equation*}
    B(\beta) =
    C_1\left(\frac{\beta}{\mu}\right)^{\frac{2}{\gamma}}\Gamma\left(1+\frac{2}{\gamma},\frac{\mu
        R_0^{\gamma}}{K\beta}\right).
  \end{equation*}
\end{collary}
\begin{proof}
  Similarly to the path loss exponent model case, we have:
  \begin{eqnarray*}
    B(\beta) &=& 2\pi\int_{0}^{\infty}rF_H((\max\{R_0,r\})^{-\gamma}(K\beta)^{-1})\d r\\
    &=&  2\pi\int_{0}^{R_0}rF_H(R_0^{-\gamma}(K\beta)^{-1})\d r +  2\pi\int_{R_0}^{\infty}rF_H(R_0^{-\gamma}(K\beta)^{-1})\d r  \\
    &=& \pi R_0^2F_H( R_0^{\gamma}(K\beta)^{-1})  + 2\pi\int_{\frac{
        R_0^{\gamma}}{\beta K}}^{\infty}p_H(t)\d t
    \int_{R_0}^{(tK\beta)^{1/\gamma}}r\d r\\
    &=&C_1\beta^{\frac{2}{\gamma}}\int_{\frac{
        R_0^{\gamma}}{K\beta}}^{\infty}t^{\frac{2}{\gamma}}p_H(t)\d t \cdotp
  \end{eqnarray*}
  We then obtain Equation (\ref{equ: corr: modif exponent}). Now
  differentiate the two sides of that equation to get:
  \begin{eqnarray*}
    B'(\beta) &=& \frac{2}{\gamma}C_1\beta^{\frac{2}{\gamma}-1}
    \int_{\frac{R_0^{\gamma}}{\beta K}}^{\infty}t^{\frac{2}{\gamma}}p_H(t)\d t+\frac{C_1R_0^2}{K^{\frac{2}{\gamma}}}p_H\left(\frac{R_0^{\gamma}}{K\beta}\right)\\
    &=& \frac{2}{\gamma}\beta^{-1}B(\beta)+\pi
    R_0^2p_H\left(\frac{R_0^{\gamma}}{K\beta}\right).
  \end{eqnarray*}
 That yields Equation (\ref{equ: corr: modif exponent path Bphay}).
  In the case of lognormal shadowing we have:
  \begin{eqnarray*}
    B(\beta) &=& C_1\beta^{\frac{2}{\gamma}}\int_{\frac{
        R_0^{\gamma}}{K\beta}}^{\infty}\frac{1}{\sqrt{2\pi\sigma_1^2}t}t^{\frac{2}{\gamma}}e^{-\frac{(\ln t)^2}{2\sigma_1^2}}\d t\\
    &=& C_1\beta^{\frac{2}{\gamma}}\int_{\ln\frac{
        R_0^{\gamma}}{K\beta}}^{\infty}\frac{1}{\sqrt{2\pi\sigma_1^2}}e^{\frac{2u}{\gamma}}e^{-\frac{u^2}{2\sigma_1^2}}\d
    u\\
    &=&C_1\beta^{\frac{2}{\gamma}}e^{(\frac{2\sigma_1}{\gamma})^2}\int_{\ln\frac{
        R_0^{\gamma}}{K\beta}}^{\infty}\frac{1}{\sqrt{2\pi\sigma_1^2}}e^{-\frac{(u-\frac{2\sigma_1^2}{\gamma})^2}{2\sigma_1^2}}\d
    u \\
    &=& C_1\beta^{\frac{2}{\gamma}}e^{(\frac{2\sigma_1}{\gamma})^2}Q\left(\frac{-\ln\beta-\ln(KR_0^{-\gamma})}{\sigma_1}-\frac{2\sigma_1}{\gamma}\right)\cdotp
  \end{eqnarray*}
  In the case of Rayleigh fading we have:
  \begin{eqnarray*}
    B(\beta) &=& C_1\beta^{\frac{2}{\gamma}}\int_{\frac{
        R_0^{\gamma}}{K\beta}}^{\infty}t^{\frac{2}{\gamma}}\mu e^{-\mu t}\d t\\
    &=&  C_1\left(\frac{\beta}{\mu}\right)^{\frac{2}{\gamma}}\Gamma\left(1+\frac{2}{\gamma},\frac{\mu
        R_0^{\gamma}}{K\beta}\right),
  \end{eqnarray*}
  by a change of  variables. Hence the results.
\end{proof}

\begin{collary}
  The number of BS covering a point $x$ is distributed according to
  the Poisson distribution of parameter $\lambda_BB(T)$. In
  particular, the outage probability given a threshold $T$ is
  \begin{equation*}
    P(\xi_0^x>T) = e^{-\lambda_BB(T)}.
  \end{equation*}
\end{collary}

\begin{proof}
  The path loss fading $Sh^x$ is a Poisson point process on
  $\R^+$ with intensity $\lambda_BB'(t)\d t$, so the number of
  point on the interval $(0,T)$ is distributed according to Poisson
  distribution of parameter $\lambda_BB(T)$.
\end{proof}

\begin{figure}
  \centering
  \includegraphics[width=300pt]{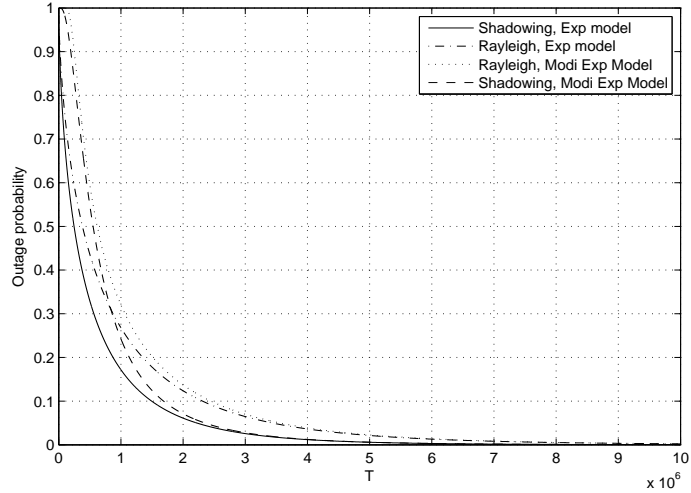}\\
  \caption{Comparison of outage probability between propagation
    models. For lognormal shadowing $\sigma = 4(dB)$, for Rayleigh
    fading $\mu=1$; $K=10^{-2}$, $\gamma=2.8$.}\label{OutageFigure}
\end{figure}

Figure \ref{OutageFigure} represents  the outage probability for
different models of fading. This  shows that the curves of modified path
loss exponent model is generally higher than those of path loss
exponent model but they are very close in the low outage region.

\section{Capacity}\label{section: Capacity user}
In this section, we calculate the mean of any capacity function  of a user. Remark that
since the system is spatially stationary the statistic of the path
loss fading and the capacity of a user does not depend on his position. Since the PDF and the CDF of the path loss fading
$\xi^x_0$ have been already calculated in Proposition \ref{proposition
  CDF PDF of a_ms}, the mean of a capacity function of a user follows
immediately. In particular:
\begin{theorem}
  The average capacity per user is
  \begin{equation}\label{equation: theorem Capacity per User}
    E(f(\xi_0^x)) =
    \lambda_B\int_0^{\infty}B'(\beta)e^{-\lambda_BB(\beta)}f(\beta)\d\beta \cdot
  \end{equation}
  In the case of path loss exponent model $L(z)=K|z|^{-\gamma}$, we
  have:
  \begin{equation}\label{equation: theorem Capacity per User expo path
      loss}
    E(f(\xi_0^x)) = \mathcal{L}_{\widetilde{f}}(\lambda_BC)
  \end{equation}
  where $\mathcal{L}_{g}(s) = \int_0^{\infty}e^{-st}g(t)\d t$ is the
  Laplace transform of the capacity function $g$ and
  $\widetilde{f}(t)=f(t^{\frac{\gamma}{2}})$.
\end{theorem}
\begin{proof}
  Equation (\ref{equation: theorem Capacity per User}) comes from
  Proposition \ref{proposition CDF PDF of a_ms}. If the path loss
  exponent model is considered, then we have:
  \begin{align*}
    E(f(\xi_0^x))
    &=\lambda_B\int_0^{\infty}\frac{2C}{\gamma}\beta^{\frac{2}{\gamma}-1}e^{-\lambda_BC\beta^{\frac{2}{\gamma}}}f(\beta)\d\beta\\
    &=\int_0^{\infty}e^{-\lambda_BC\beta_1}f(\beta_1^{\frac{\gamma}{2}})\d\beta_1
  \end{align*}
  by the change of  variable $\beta_1 = \beta^{\frac{\gamma}{2}}$.
\end{proof}
The statistic of the cell capacity $S_o(f)$ is more difficult to
analyze. In this section, we calculate its mean $m(f)$ and  lower
bound and upper bound of its variance $v(f)$. We state the following
lemma, which is straightforward due to Assumption \ref{assumption} but
still useful:
\begin{lemma}\label{lemma: bex bey indepent}
  Given a fixed configuration $\Pi_B$ of BSs, the Poisson point
  processes of path loss fading $Sh^x$ and $Sh^y$ are independent
  for any two different points $x,\, y$.
\end{lemma}
\begin{lemma}\label{lemma: pdf of P_yx}
Let $z=y-x$. The PDF of $s_{yx}$ is given
  by
  \begin{equation*}
    p_{s_{yx}}(t) = \frac{1}{l(z)t^2}p_H\left(\frac{1}{L(z)t}\right).
  \end{equation*}
\end{lemma}
\begin{proof}
  We have
  \begin{eqnarray*}
    P(s_{yx}<t)   &=& P(h_{yx}>\frac{1}{L(z)t})\\
    &=&  F_H\left(\frac{1}{L(z)t}\right).
  \end{eqnarray*}
  The density probability function is then
  \begin{eqnarray*}
    p_{s_{yx}}(t) = \frac{1}{L(z)t^2}p_H\left(\frac{1}{L(z)t}\right).
  \end{eqnarray*}
\end{proof}
\begin{theorem}\label{theorem: expectation of Sf}
  The expectation of the cell capacity of the typical BS is
  \begin{equation}\label{equation: expectation of Sf}
    m(f) =
    \lambda_M\int_0^{\infty}B'(\beta)e^{-\lambda_BB(\beta)}f(\beta)\d\beta \cdot
  \end{equation}
  In the case of path loss exponent model $L(z)=K|z|^{-\gamma}$, we
  have:
  \begin{equation}\label{equation: theorem sum Capa Exp Pathloss}
    m(f) =
    \frac{\lambda_M}{\lambda_B}\mathcal{L}_{\widetilde{f}}(\lambda_BC)\cdot
  \end{equation}
\end{theorem}
\begin{proof}
  Given a fixed configuration of BSs $\Pi_B$, the random variables
  $1(s_{ox}<\xi_0^x)f(s_{ox})$ obtained from all $x\in \mathbf{R}^2$
  are independent. Thus, the marked point process
  $\tilde{\Pi}_M=(x_i,1(s_{ox_i}<\xi_0^{x_i})f(s_{ox_i}))$ is a
  Poisson point process. Using the Campbell theorem we have:
  \begin{equation*}
    E(S_o(f)\mid \Pi_B) =
    \lambda_M\int_{\mathbf{R}^2}E\left(1(s_{ox}<\xi_0^x)f(s_{ox})\mid \Pi_B\right)\d x \cdot
  \end{equation*}
  As a consequence,
  \begin{align*}
    E(S_o(f)) &= E\left(\lambda_M\int_{\mathbf{R}^2}E(1(s_{ox}<\xi_0^x)f(s_{ox})\mid \Pi_B)\d x\right)\\
    &=\lambda_M\int_{\mathbf{R}^2}E\left(1(s_{ox}<\xi_0^x)f(s_{ox})\right)\d x \cdot
  \end{align*}
  In virtue of  Lemma \ref{lemma: pdf of P_yx}, Proposition
  \ref{proposition: POISSON POINT PROCESS of path shadow} and Collary
  \ref{proposition CDF PDF of a_ms}, we have:
  \begin{align*}
    E(S_o(f)) &= \lambda_M\int_{\mathbf{R}^2}E(1(s_{ox}<\xi_0^x)f(s_{ox}))\d x\\
    &=
    \lambda_M\int_{\mathbf{R}^2}\int_0^{\infty}p_{s_{ox}}(t)P(t<\xi_0^x)f(t)\d t\d x\\
    &=
    \lambda_M\int_0^{\infty}f(t)e^{-\lambda_BB(t)}\d t\int_{\mathbf{R}^2}\frac{\partial
      F_H}{\partial t}\left(\frac{1}{L(x)t}\right)\d x\\
    &= \lambda_M\int_0^{\infty}f(t)e^{-\lambda_BB(t)}\d t
    \frac{\partial
    }{\partial t}\left(\int_{R^2}F(\frac{1}{L(x)t})\d x\right)\\
    &=\lambda_M\int_0^{\infty}f(t)e^{-\lambda_BB(t)}B'(t)\d t    \cdot
  \end{align*}
  For the case of path loss exponent model, Equation \ref{equation:
    theorem sum Capa Exp Pathloss} follows easily. This completes the
  proof.
\end{proof}
Equation (\ref{equation: expectation of Sf}) has the following
interpretation: the mean cell capacity  is the product of the
mean number of users per cell and the mean  capacity per user.
\begin{theorem} \label{theorem: lower bound covariance S_o(f) S_o(g)}
  Given two capacity functions $f,g$ we have~:
  \begin{equation}\label{equation: theorem : covariance f g > mean fg
    }
    \cov(S_o(f),S_o(g))\geq m(f.g)   \cdot
  \end{equation}
In particular,
\begin{equation*}
  \var(S_o(f))\geq m(f^2).
\end{equation*}
\end{theorem}
\begin{proof}
  For simplicity, let $\beta_x = s_{ox}1(s_{ox}<\xi_0^x)$ and $f(0)=0,
  g(0)=0$, we have :
  \begin{align*}
    \cov(S_o(f)S_o(g)) &= E(\cov(S_o(f),S_o(g)\mid \Pi_B))  \\
    & + E(E(S_o(f)\mid \Pi_B)E(S_o(g)\mid \Pi_B)) - E(S_o(f))E(S_o(g))\\
    &= T_1 + T_2 - T_3 \cdot
  \end{align*}
  It is clear that
  \begin{eqnarray*}
    T_3 = m(f)m(g) \cdot
  \end{eqnarray*}
  Consider the first term. Remind that we have assumed that all random
  fading $\{h_{yx}\}_{y,x\in \mathbf{R}^2}$ are independent, so given
  a fixed configuration $\Pi_B$ of BSs, the random variables
  $\{\beta_x\}_{x\in \mathbf{R}^2}$ are independent. Hence by Campbell
  formula we have:
  \begin{eqnarray*}
    T_1 &=& \lambda_ME\int_{\mathbf{R}^2}E(f(\beta_x)g(\beta_x)\mid \Pi_B)\d x\\
    &=&  \lambda_M\int_{\mathbf{R}^2}E(E(f(\beta_x)g(\beta_x)\mid \Pi_B))\d x\\
    &=&  \lambda_M\int_{\mathbf{R}^2}E(f(\beta_x)g(\beta_x))\d x\\
    &=&  m(f.g)   \cdot
  \end{eqnarray*}
  Now consider the second term
  \begin{multline*}
    T_2 = \lambda_M^2E\left(\int_{\mathbf{R}^2}E(f(\beta_x)\mid \Pi_B)\d x\int_{R^2}E(g(\beta_x)\mid \Pi_B)\d x\right)\\
    \begin{aligned}
      &= \lambda_M^2E\left(\int_{\mathbf{R}^2}\int_{\mathbf{R}^2}E(f(\beta_{x})\mid \Pi_B)E(g(\beta_{y})\mid \Pi_B)\d x\d y\right)\\
      &= \lambda_M^2\int_{\mathbf{R}^2}\int_{\mathbf{R}^2}E\left(E(f(\beta_{x})\mid \Pi_B)E(g(\beta_{y})\mid \Pi_B)\right)\d x\d y\\
      &=\lambda_M^2\int_{\mathbf{R}^2}\int_{\mathbf{R}^2}E(f(\beta_{x}))E(g(\beta_{y}))\d x\d y\\
      &=\lambda_M^2\int_{\mathbf{R}^2}\int_{\mathbf{R}^2}\int_0^{\infty}\int_0^{\infty}P(s_{ox}<\xi_0^x,
      s_{oy}<\xi_0^y\mid s_{ox}=t_1,s_{oy}=t_2)\times
    \end{aligned}\\
    \times f(t_1)g(t_2)p_{s_{ox}}(t_1)p_{s_{oy}}(t_2) \d t_1\d t_2\d x\d y
    \cdot
  \end{multline*}
  by remarking that $\beta_{x}$ and $\beta_{y}$ are independent if
  $x\neq y$ (Lemma \ref{lemma: bex bey indepent}). We will prove that
  if $x\neq y$:
  \begin{multline*}
    P(s_{ox}<\xi_0^x, s_{oy}<\xi_0^y\mid s_{ox}=t_1,s_{oy}=t_2)\geq\\
    P(s_{ox}<\xi_0^x\mid s_{ox}=t_1)P(s_{oy}<\xi_0^y\mid s_{oy}=t_2)
  \end{multline*}
  Consider the marked point process
  $\Pi_B^{x,y}=\{y_i,h_{y_ix},h_{y_iy})\}$. Since the marks are
  independent, it is a Poisson point process on $\mathbf{R}^4$ with
  intensity
$$m_a^{x,u_1,u_2}(\d y,du_1,du_2)=\lambda_M\d y\otimes
p_H(u_1)du_1\otimes p_H(u_2)du_2\cdot$$ Consider two sets
\[A_1 = \{(y,u_1,u_2): L(y)u_1\geq t_1^{-1}\}\] and
\[A_2= \{(y,u_1,u_2): L(y)u_2\geq t_2^{-1}\},\] we have:
\begin{multline*}
  P(s_{ox}<\xi_0^x, s_{oy}<\xi_0^y\mid s_{ox}=t_1,s_{oy}=t_2)
  =  P(\Pi_B^{x,y}(A_1\cup A_2)=\emptyset) \\
  \begin{aligned}
    &   =  e^{-m_a^{x,y}(A_1\cup A_2)} \\
    &   \geq e^{-m_a^{x,y}(A_1)-m_a^{x,y}(A_2)}\\
    &  = P(\Pi_B^{x,y}( A_1)=\emptyset)P(\Pi_B^{x,y}( A_2)=\emptyset) \\
    & = P(s_{ox}<\xi_0^x\mid s_{ox}=t_1)P(s_{oy}<\xi_0^y\mid
    s_{oy}=t_2)\cdot
  \end{aligned}
\end{multline*}
Thus,
\begin{align*}
  T_2 &\geq
  \lambda_M^2\int_{\mathbf{R}^2}\int_{\mathbf{R}^2}\int_0^{\infty}\int_0^{\infty}P(s_{ox}<\xi_0^x\mid
  s_{ox}=t_1)P(s_{oy}<\xi_0^y\mid
  s_{oy}=t_2)\times\\
  & \ \ \ \ \times f(t_1)g(t_2)p_{s_{ox}}(t_1)p_{s_{oy}}(t_2) \d t_1\d t_2 \d x\d y\\
  &= m(f)m(g)\cdot
\end{align*}
The result follows.
\end{proof}
\begin{theorem}\label{theorem: Upper bound on the covariance of sum} \label{proposition: Upper bound on the variance of
    sum} For $f$ and $g$ two capacity functions, we have:
  \begin{equation*}
    \cov(S_o(f),S_o(g))\leq m(f.g) + m(f)n(g) -m(f)m(g)
  \end{equation*}
  where
  \begin{equation}\label{equation: expectation n(f)}
    n(f) =
    \lambda_M\int_0^{\infty}B'(t)f(t)\d t  \cdot 
  \end{equation}
\end{theorem}
\begin{proof}
  We continue the proof of Theorem \ref{theorem: lower bound
    covariance S_o(f) S_o(g)}, we have to prove that $$T_2\leq
  m(f)n(g)\cdot$$ Indeed,
  \begin{eqnarray*}
    P(s_{ox}<\xi_0^x, s_{oy}<\xi_0^y\mid
    s_{ox}=t_1,s_{oy}=t_2)\leq P(s_{ox}<\xi_0^x\mid s_{ox}=t_1),
  \end{eqnarray*}
  thus,
  \begin{align*}
    T_2
    &\leq\lambda_M^2\int_{\mathbf{R}^2}\int_{\mathbf{R}^2}\int_0^{\infty}\int_0^{\infty}P(s_{ox}<\xi_0^x\mid
    s_{ox}=t_1)\times\\
    &\ \ \ \ \ \times f(t_1)g(t_2)p_{s_{ox}}(t_1)p_{s_{oy}}(t_2) \d t_1\d t_2 \d x\d y\\
    &\leq m(f)\int_{0}^{\infty}g(t_2)\int_{\R^2}\frac{\partial
      F_H}{\partial t}(\frac{1}{L(x)t}) \d x \d t_2\\
    &= m(f)n(g) \cdot
  \end{align*}
  Hence the result.
\end{proof}
\section{Examples} \label{section: Examples} 
\subsection{Number of users in a cell}\label{subsection: number of
  users per cell}
For $f_0(t)=1$, the random variable
$n_o:=S(f_0)=\sum_{i=0}^{\infty}1(x_i\in C_o)$ represents the number
of users who view $o$ as the best server, and thus will be served by $o$.
\begin{eqnarray*}
  E(n_o) &=&
  \lambda_M\int_0^{\infty}B'(\beta)e^{-\lambda_BB(\beta)}d\beta\\
  &=& \frac{\lambda_M}{\lambda_B} \cdot
\end{eqnarray*}
The mean number of users served by a BS is
$\frac{\lambda_M}{\lambda_B}$ which is easily interpreted. We rewrite
the formula \eqref{equation: expectation of Sf} by
\begin{equation*}
  E(S_o(f)) = \frac{\lambda_M}{\lambda_B}E(f(\xi_{0}^x)).
\end{equation*}
Again this is easily interpreted. The average sum rate is the product
of the average user per cell and the average per user.

Now apply Theorem  \ref{theorem: lower bound covariance S_o(f) S_o(g)}, we get that
\begin{equation*}
  \var(n_o) \geq m(1) = \frac{\lambda_M}{\lambda_B}\cdotp
\end{equation*}
We can not apply Theorem \ref{proposition: Upper bound on the
  variance of sum} because
$n(f_0)=\lambda_M\int_0^{\infty}B'(t)\d t=\infty$.
\begin{figure}
  \centering
  \includegraphics[width=300pt]{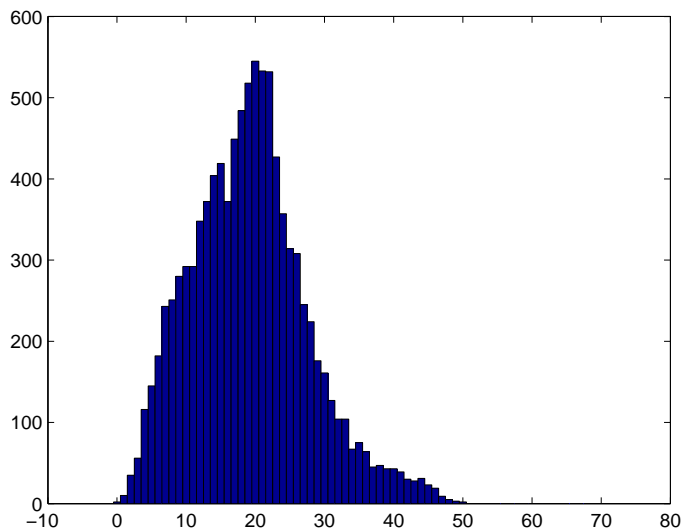}\\
  \caption{Histogram of $n_o$} \label{HistNoFigure}
\end{figure}

\subsection{Number of users in outage in a cell}
Consider $f_1(t)=1(t> T)$, then $S_o(f_1)$ is the number of users in
outage in the typical cell. We have
\begin{align*}
  m(f_1) &=
  \lambda_M\int_T^{\infty}B'(\beta)e^{-\lambda_BB(\beta)}d\beta\\
  &=\frac{\lambda_M}{\lambda_B}e^{-\lambda_BB(T)},
\end{align*}
and
\begin{eqnarray*}
  v(f_1) \geq \frac{\lambda_M}{\lambda_B}e^{-\lambda_BB(T)} \cdot
\end{eqnarray*}
Note that again,  we can not apply Theorem~\ref{proposition:
  Upper bound on the variance of sum} as
$n(f_1)=\lambda_M\int_T^{\infty}B'(t)\d t$ is infinite.

\subsection{Number of covered users in a cell}
Consider $f_2(t)=1(t\leq T)$, then $S(f_2)$ represents the number of
covered users in the typical cell. We have:
\begin{align*}
  m(f_2) &=
  \lambda_M\int_0^{T}B'(\beta)e^{-\lambda_BB(\beta)}d\beta\\
  &=\frac{\lambda_M}{\lambda_B}\left(1-e^{-\lambda_BB(T)}\right),
\end{align*}
\begin{eqnarray*}
  v(f_2) \geq \frac{\lambda_M}{\lambda_B}\left(1-e^{-\lambda_BB(T)}\right),
\end{eqnarray*}
and
\begin{equation*}
  v(f_2) \leq \frac{\lambda_M}{\lambda_B}e^{-\lambda_BB(T)}+\frac{\lambda_M^2}{\lambda_B^2}e^{-\lambda_BB(T)}\left(\lambda_BB(T)-1+e^{-\lambda_BB(T)}\right)\cdot
\end{equation*}

\subsection{Total bit rate of a cell}
We now consider the piecewise constant function $f_3(t) =
\sum_1^n1(T_i\leq t<T_{i+1})c_i$ with $0<T_1<T_2<...<T_n<T_{n+1}$ and
$T_{n+1}$ can be infinite. If $f_3$ is the function that represents
the actual bit rate then $S_o(f_3)$ represents the total bit rates of
all users in the cell. We have:

\begin{eqnarray*}
  m(f_3) &=&
  \lambda_M\int_0^{\infty}B'(\beta)e^{-\lambda_BB(\beta)}\sum_1^n1(T_i\leq \beta<T_{i+1})c_id\beta\\
  &=& \lambda_M \sum_{i=1}^nc_i\left( e^{-\lambda_BB(T_i)}-e^{-\lambda_BB(T_{i+1})}\right) \cdot
\end{eqnarray*}

\begin{eqnarray*}
  v(f_3) &\geq& \lambda_M \sum_{i=1}^nc_i^2\left( e^{-\lambda_BB(T_i)}-e^{-\lambda_BB(T_{i+1})}\right) \cdot
\end{eqnarray*}

\begin{eqnarray*}
  v(f_3) &\leq& \lambda_M \sum_{i=1}^nc_i^2\left( e^{-\lambda_BB(T_i)}-e^{-\lambda_BB(T_{i+1})}\right) +\\
  &&  + \left(\frac{\lambda_M}{\lambda_B}\right)^2\sum_{i=1}^n\left( e^{-\lambda_BB(T_i)}-e^{-\lambda_BB(T_{i+1})}\right)\times\\
  \\&& \times \sum_{i=1}^n\left( \lambda_BB(T_{i+1})-\lambda_BB(T_{i})- e^{-\lambda_BB(T_i)}+e^{-\lambda_BB(T_{i+1})}\right) \cdot
\end{eqnarray*}
\subsection{Discussion on the distribution of $S_o(f)$}
\begin{figure}
  \centering
  \includegraphics[width=300pt]{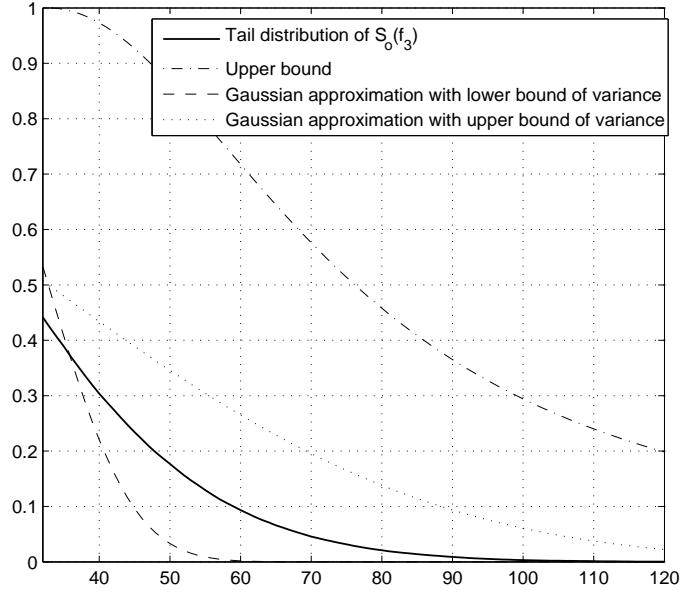}\\
  \caption{Tail distribution of $S_o(f_3)$} \label{Tail Distribution}
\end{figure}
\begin{figure}
  \centering
  \includegraphics[width=300pt]{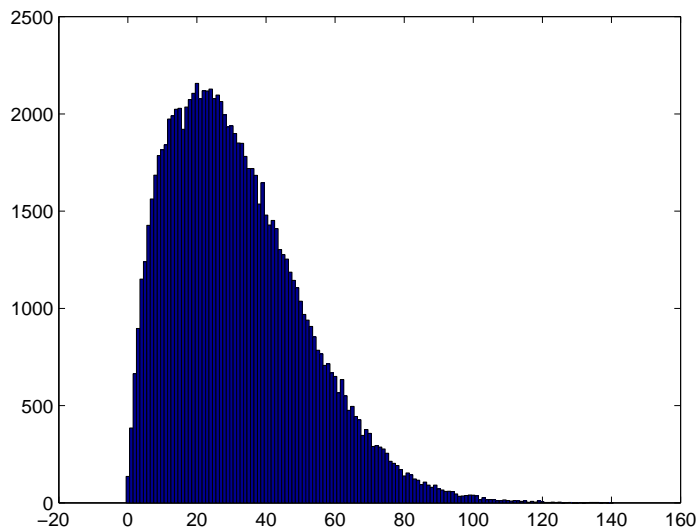}\\
  \caption{A typical histogram of $S_o(f_3)$} \label{Histogram Sf3}
\end{figure}

The distribution of $S_o(f)$ does not behave like a Gaussian
distribution even in the limit regimes. Take, for example, the
histogram of $n_o=S_o(f_0)$ and that of $S_o(f_3)$ which are shown in
figures \ref{HistNoFigure} and \ref{Histogram Sf3} respectively. For
the case of no fading, $H=$ constant, in
\cite{OnacertainVoronoiaggregativeprocessrelatedtoabivariatePoissonprocess}
the author found some approximative but not reliable bounds of the
distribution of $S_o(f)$ for equivariant functions $f$ but no
approximation or bounds is found for general capacity functions. In
addition, no closed expression is found for the Laplace transform of
functional $S_o(f)$. In our case where the fading is considered, this
is expected to be more challenging.

We can find an upper bound for the tail distribution by Chebyshev's
inequality:
\begin{eqnarray*}
  P(S_o(f)>m(f)+t)&\leq& \frac{v(f)}{v(f)+t^2}\\
  &\leq& \frac{m(f^2) + m(f)n(f) - (m(f))^2}{m(f^2) + m(f)n(f) -
    (m(f))^2+t^2}
\end{eqnarray*}
The above inequality provides a robust upper bound for the tail
distribution and valid for all capacity function $f$. However the gap
is large (Figure \ref{Tail Distribution}). It is well known that other
types of concentration inequality based on Chernoff bound can give
better bound. In this direction,
\cite{ConcentrationAndDeviationInequalitiesInInfiniteDimensionsViaCovarianceRepresentations},
\cite{UpperBoundsOnRubinsteinDistancesOnConfigurationSpaces} and
\cite{ANewModifiedLogarithmicSobolevInequalityForPoissonPointProcessesAndSeveralApplications}
provide concentration inequalities that apply for functional related
to one PPP. These inequalities can not be directly applied in our case
because our target is a functional related to two independent
PPPs. Actually we can combine the two independent PPPs into one united
PPP by the independent marking theorem. Unfortunately the functional
$S_o(f)$ of the united PPP does not satisfy the required conditions
for the concentration inequalities neither on
\cite{ConcentrationAndDeviationInequalitiesInInfiniteDimensionsViaCovarianceRepresentations},
\cite{ANewModifiedLogarithmicSobolevInequalityForPoissonPointProcessesAndSeveralApplications}
nor on
\cite{UpperBoundsOnRubinsteinDistancesOnConfigurationSpaces}. But we
believe that similar techniques used in these references can be used
to derive a upper bound the tail distribution of $S_o(f)$.
\section{Conclusion} In this paper we introduce a general model to
evaluate the outage probability and the capacity of wireless noise
limited network. It is in fact an extension of models introduced in
series of papers
\cite{Stochasticgeometryandarchitectureofcommunicationnetworks},
\cite{StochasticgeometryandarchitectureofcommunicationnetworksFrench},
\cite{OnacertainVoronoiaggregativeprocessrelatedtoabivariatePoissonprocess}. The
main difference is that we take into account the effect of fading, and
that we assume that a user connects to the BS with strongest signal
rather than the closest one. We first show that for a particular
user, the path loss fading process from all BSs seen from this user is
a Poisson point process in the positive half line. We find
explicit expression for the outage probability, the expectation of
capacity of a user, and the expectation of the cell capacity of the
typical BS $S_o(f)$. We find the lower bound and upper bound for the
variance of the cell capacity. We consider general model for path loss
and fading. The results presented in this paper actually generalizes
the results on
\cite{OnacertainVoronoiaggregativeprocessrelatedtoabivariatePoissonprocess}. Possible
further research is to find a way to compute the distribution of
$S_o(f)$.







\end{document}